\documentclass[11pt]{amsart}
\usepackage{amsmath, amsthm, latexsym, amssymb, amsfonts, epsfig, epsf}

\newenvironment{pf}{\proof[\proofname]}{\endproof}
\theoremstyle{plain}
\newtheorem*{Th1}{Theorem 1}
\newtheorem*{Th2}{Theorem 2}
\newtheorem{Th}{Theorem}[section]

\newtheorem{Prop}[Th]{Proposition}
\newtheorem{Lemma}[Th]{Lemma}
\numberwithin{equation}{section} \theoremstyle{definition}
\newtheorem{Rem}[Th]{Remark}

\newtheorem{Ex}{Example}

\newtheorem{Def}[Th]{Definition}

\newcommand{\cal}[1]{\mathcal{#1}}
\newcommand{\F}{\mathbb F}

\newcommand{\Z}{\mathbb Z}
\newcommand{\R}{\mathbb R}

\newcommand{\D}{\Delta}

\newcommand{\cC}{\cal C}

\newcommand{\cL}{\cal L}

\newcommand{\Sig}{\Sigma}

\newcommand{\spn}{\operatorname{span}}

\newcommand{\GL}{\operatorname{GL}}

\newcommand{\A}[1]{\operatorname{AGL}(#1,\Z)}

\newcommand{\rs}[1]{Section~\ref{S:#1}}
\newcommand{\rl}[1]{Lemma~\ref{L:#1}}
\newcommand{\rp}[1]{Proposition~\ref{P:#1}}

\newcommand{\rr}[1]{Remark~\ref{R:#1}}
\newcommand{\rex}[1]{Example~\ref{ex:#1}}
\newcommand{\re}[1]{(\ref{e:#1})}

\newcommand{\rt}[1] {Theorem~\ref{T:#1}}
\newcommand{\rd}[1]{Definition~\ref{D:#1}}
\newcommand{\rf}[1]{Figure~\ref{F:#1}}

\begin{document}

\author{ Ivan Soprunov, Jenya Soprunova}
\title[Toric surface codes and Minkowski length of polygons]
{Toric surface codes and\\ Minkowski length of polygons}

\dedicatory{To our advisor Askold Khovanskii on the occasion of his 60-th anniversary, with love.}

\begin{abstract}
In this paper we prove new lower bounds for the minimum
distance of a toric surface code $\cC_P$ defined by a convex lattice
polygon $P\subset \R^2$.
The bounds involve a geometric invariant $L(P)$, called the full Minkowski
length of $P$ which
can be easily computed for any given $P$.
\end{abstract}

\maketitle




\section*{Introduction}

Consider a convex polygon $P$ in $\R^2$ whose vertices lie in the
integer lattice $\Z^2$.
It determines a vector space $\cL_K(P)$ (over a filed $K$) of
polynomials $f(t_1,t_2)$ whose monomials
correspond to the lattice points in $P$:
$$\cL_K(P)=\spn_K\{t_1^{m_1}t_2^{m_2}\ |\ (m_1,m_2)\in P\cap\Z^n\}.$$
Consider a finite field $\F_q$.
The {\it toric surface code} $\cC_P$, first introduced by Hansen in
\cite{Han}, is defined
by evaluating the polynomials in $\cL_{\F_q}(P)$ at all the points
$(t_1,t_2)$ in the algebraic torus
$(\F_q^*)^2$. To be more precise, $\cC_P$ is a linear code whose codewords 
are the strings $(f(t_1,t_2)\ |\ (t_1,t_2)\in (\F_q^*)^2)$ for
 $f\in\cL_{\F_q}(P)$.
 It is convenient to assume that $P$ is contained in the square $K^2_q=[0,q-2]^2$
so that all the monomials in $\cL_{\F_q}(P)$ are linearly independent
over $\F_q$. Thus $\cC_P$ has block length  $(q-1)^2$ and 
dimension equal to the number of the lattice points in $P$.

Note that the weight of each non-zero codeword in $\cC_P$ is
the number of points $(t_1,t_2)\in(\F_q^*)^2$ where the corresponding polynomial does not
vanish. Therefore, the minimum distance of $\cC_P$ (which is the minimum weight for linear codes)
equals 
$$d(\cC_P)=(q-1)^2-\max_{0\neq f\in\cL_{\F_q}(P)}Z(f),$$
where $Z(f)$ is the number of zeroes  (i.e. points of vanishing) in $(\F_q^*)^2$ 
of $f$.

The name {\it toric surface code} comes from the
fact that $P$ defines a toric surface $X$ over $\overline\F_q$
(strictly speaking the fan that
defines $X$ is a refinement of the normal fan of~$P$), where
$\cL_{\overline\F_q}(P)$ can be identified
with the space of global sections of a semiample divisor on $X$ (see
for example \cite{Fu}).
This allows to exploit algebraic geometric techniques to produce
results about the minimum distance
of $\cC_P$. In particular, Little and Schenck in \cite{LiSche} used
intersection theory on
toric surfaces to come up with the following general idea: If $q$ is
sufficiently large then polynomials
$f\in\cL_{\F_q}(P)$ with more absolutely irreducible factors will
necessarily have more
zeroes in $(\F_q^*)^2$ (\cite{LiSche}, Proposition 5.2).

In this paper we expand this idea to produce explicit bounds for the
minimum distance
of $\cC_P$ in terms of certain geometric invariant $L(P)$, which we
call the {full Minkowski
length} of $P$. Essentially $L(P)$ tells you the largest possible number of
absolutely irreducible factors a polynomial $f\in\cL_{\F_q}(P)$ can have, but it
derives it from the geometry of the polygon $P$ (see \rd{maximal}).
The number $L(P)$
is easily computable --- we give a simple algorithm which is
polynomial in the number of lattice points in $P$. Moreover we obtain
a description of the factorization $f=f_1\cdots f_{L(P)}$
for $f\in\cL_{\F_q}(P)$ with the largest number of factors. More
precisely, in  \rp{HWL} we show
that the Newton polygon $P(f_i)$ (which is the convex hull of the
exponents of the
monomials in $f_i$) is either a primitive segment, a unit simplex, or
a triangle with exactly 1 interior
and 3 boundary lattice points, called an {\it exceptional triangle}.
This description enables us to  prove the following bound:

\begin{Th1} Let $P\subset K_{q}^2$ be a lattice polygon with area $A$
and full Minkowski length~$L$. Then for $q\geq
\max\left(23,\big(c+\sqrt{c^2+5/2}\big)^2\right)$, where $c=A/2-L+9/4$,
the minimum distance of the toric surface code $\cC_P$ satisfies
$$d(\cC_P)\geq (q-1)^2-L(q-1)-2\sqrt{q}+1.$$
\end{Th1}

The condition that no factorization $f=f_1\cdots f_{L(P)}$  contains
an exceptional triangle (as the Newton polygon of one of the factors) is
geometric and can be easily checked for any given $P$
(we provide a simple algorithm for this which is polynomial in the
number of lattice points in $P$).
In this case we have a better bound for the minimum distance of the
toric surface code:

\begin{Th2} Let $P\subset K_{q}^2$ be a lattice polygon with area $A$
and full Minkowski
length $L$. Under the above condition on $P$,
for $q\geq \max\left(37,\big(c+\sqrt{c^2+2}\big)^2\right)$, where
$c=A/2-L+11/4$,
the minimum distance of the toric surface code $\cC_P$ satisfies
$$d(\cC_P)\geq (q-1)^2-L(q-1).$$
\end{Th2}

We remark that our thresholds for $q$ where the bounds begin to hold
are much smaller than the
ones in Little and Schenck's result (\cite{LiSche}, Proposition 5.2).

Although, as mentioned above, the minimum distance problem for toric
codes is tightly connected to
toric varieties, all our methods are geometric and combinatorial and
do not use algebraic geometry (except for the Hasse--Weil bound, see
\rs{prelim}).
In \rs{length} we define the full Minkowski length $L(P)$ and
establish combinatorial
properties of polygons with $L(P)=1,2$. In \rs{main} we give a proof
of Theorem 1 and Theorem 2.
\rs{alg} is devoted to the above mentioned algorithms for computing
$L(P)$ and determining the presence of an exceptional triangle.
Finally, in
\rs{ex} we give a detailed analysis of three toric surface codes which
illustrates our methods.

\subsection*{Acknowledgments} We thank Leah Gold and Felipe Martins
for helpful discussions
on coding theory.


\section{Full Minkowski length of polytopes}\label{S:length}


\subsection{Minkowski sum}
Let $P$ and $Q$ be convex polytopes in $\R^n$. Their {\it Minkowski sum} is
$$P+Q=\{p+q\in\R^n\ |\ p\in P,\ q\in Q\},$$
which is again a convex polytope. \rf{sum} shows the Minkowski sum of a triangle and a square.
\begin{figure}[h]
\centerline{
 \scalebox{0.6}
 {
\input{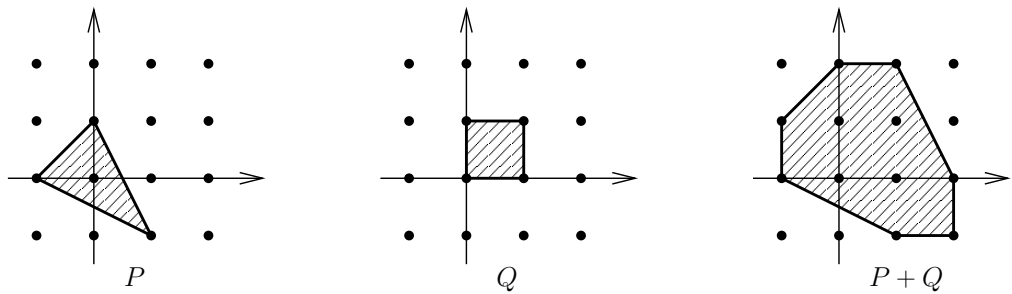}}}
\caption{The Minkowski sum of two polygons}
\label{F:sum}
\end{figure}

Let $f$ be a Laurent polynomial in  $K[t_1^{\pm 1},\dots,t_n^{\pm 1}]$ (for some field $K$). Then its
{\it Newton polytope} $P_f$ is the convex hull of the exponent vectors 
of the monomials appearing in $f$. Thus $P_f$ is a {\it lattice polytope} as its vertices belong to the integer lattice $\Z^n\subset\R^n$. Note that if $f,g\in K[t_1^{\pm 1},\dots,t_n^{\pm 1}]$ then the Newton polytope of their product $P_{fg}$ is the Minkowski sum $P_f+P_g$.  A {\it primitive lattice segment} $E$
is a line segment whose only lattice points are its endpoints. The difference of the endpoints 
is a vector $v_E$ whose coordinates are relatively prime ($v_E$ is defined up to sign).
A polytope which is the Minkowski sum of primitive lattice segments is called a {\it (lattice) zonotope}. 

The automorphism group of the lattice is the group of {affine unimodular transformations}, denoted 
by  $\A{n}$, which consists of translations by an integer vector and linear transformations in $\GL(n,\Z)$. Affine unimodular transformations correspond to monomial changes of variables in $K[t_1^{\pm 1},\dots,t_n^{\pm 1}]$ and preserve the zero set of $f$ in the algebraic torus $(K^*)^n$. 

\subsection{Full Minkowski length}
Let $P$ be a lattice polytope in $\R^n$.  Consider a Minkowski decomposition
$$P=P_1+\dots+P_\ell$$
into lattice polytopes $P_i$ of positive dimension. Clearly, there are
only finitely many such decompositions.  We let $\ell(P)$ be the largest number of
summands in such decompositions of $P$, and call it the {\it Minkowski length}
of $P$.

\begin{Def}\label{D:maximal}
The {\it full Minkowski length} of $P$ is the maximum of the
 Minkowski lengths of all subpolytopes $Q$ in $P$,
 $$L(P):=\max\{\ell(Q)\, | \, Q\subseteq P\}.$$
A subpolytope $Q\subseteq P$ is called {\it maximal} for $P$ if
$\ell(Q)=L(P)$. A Minkowski decomposition of $Q$ into $L(P)$ summands
of positive dimension will be referred to as  a {\it maximal
(Minkowski) decomposition in P}.
\end{Def}

Here are a few simple properties of $L(P)$ and maximal subpolytopes.

\begin{Prop}\label{P:simple}  Let $P$, $P_1$, $P_2$, and $Q$ be
lattice polytopes in $\R^n$.
\begin{enumerate}
\item $L(P)$ is $\A{n}$-invariant.
\item $L(P)\geq 1$ if and only if $\dim(P)>0$.
\item If $P_1+P_2\subseteq P$ then $L(P_1)+L(P_2)\leq L(P)$.
\item If $Q$ is maximal for $P$ then $Q$ contains a zonotope $Z$
maximal for $P$.
\end{enumerate}
\end{Prop}

\begin{pf} The first three statements are trivial. For the forth one,
note that if
$$Q=Q_1+\dots+Q_{L(P)}$$
is a maximal Minkowski decomposition in $P$
then by replacing each $Q_i$ with one of its edges we obtain  a
zonotope $Z\subseteq Q$ with $\ell(Z)\geq L(P)$. But $Z\subseteq P$,
so $\ell(Z)=L(P)$.
\end{pf}

Notice that the summands of every maximal decomposition in $P$ are polytopes of
full Minkowski length 1. It seems  to be a hard problem to describe
polytopes of full Minkowski length 1  in general. However, in
dimensions 1 and  2 we do have a simple description for such
polytopes (\rt{indecomp}).

\begin{Def} A lattice polytope $P$ is {\it strongly indecomposable} if
its full Minkowski length $L(P)$ is 1.
In other words, no subpolytope $Q\subseteq P$ is a Minkowski sum of
lattice polytopes of positive dimensions.
 \end{Def}

Clearly, primitive segments are strongly indecomposable and they are
the only 1-dimensional strongly indecomposable polytopes.

 Let $\D$ be the standard 2-simplex and $T_0$ be the triangle with
vertices $(1,0)$, $(0,1)$ and $(3,3)$ (see \rf{indecomp}).  It is easy
to see that the they both are strongly indecomposable.


 \begin{figure}[h]
\centerline{
 \scalebox{0.6}
 {
\input{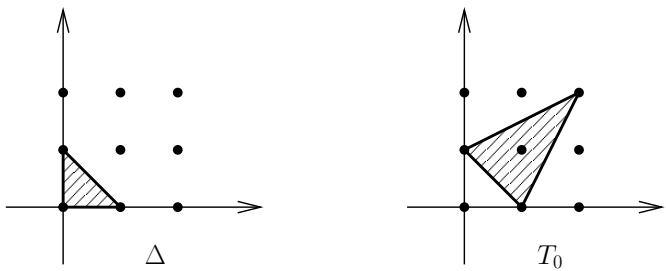}}}
\caption{Strongly indecomposable polygons}
\label{F:indecomp}
\end{figure}

 Next theorem shows that these are essentially the only strongly
indecomposable polygons. In the proof of this theorem and
frequently later in the paper we will use Pick's formula: Let $P$ be
a lattice polygon in $\R^2$. Then the area of $P$ equals
$$A=I+\frac{B}{2}-1,$$
where $I$ is the number of interior lattice points in $P$ and $B$ is the
number of boundary points in $P$. The proof of this formula  can
be found for example in \cite{BeRo}.

\begin{Th}\label{T:indecomp}
 Let $P$ be a strongly indecomposable polygon. Then $P$ is $\A2$-equivalent
 to either the standard 2-simplex $\D$ or the triangle $T_0$ above.
 \end{Th}

\begin{pf}
First, note that $P$ cannot contain more than 4 lattice points.
Indeed, suppose $a=(a_1,a_2)$
and $b=(b_1,b_2)$ lie in $P\cap\Z^2$. If $a_i\equiv b_i \mod 2$, for $i=1,2$,
then the segment $[a,b]$ lies in $P$ and is not primitive, hence,
$L(P)>1$. Since
there are only 4 possible pairs of remainders mod $2$, $P$ has at most
4 lattice points.

Suppose $P$ is a triangle, then its  sides must be primitive and
either $P$ has no interior lattice points or it has exactly one
interior lattice point. In the first case $P$ has area $1/2$ (by
Pick's formula) and so is $\A{2}$-equivalent to $\D$. In the second
case $P$ has area $3/2$ (by Pick's formula) and hence any two of its
sides generate a parallelogram of area 3.
Every such triangle is  $\A{2}$-equivalent to~$T_0$.

Now suppose $P$ is a quadrilateral. Then it has no interior lattice
points and so its area is 1 (by Pick's formula). Every such
quadrilateral is $\A{2}$-equivalent to the unit square. But the unit
square is  obviously decomposable.
\end{pf}

\begin{Def} A lattice polygon is called a  {\it unit triangle} if it
is $\A{2}$-equivalent to $\D$, and an {\it exceptional triangle} if it is
 $\A{2}$-equivalent to $T_0$.
\end{Def}

The following theorem describes maximal Minkowski decompositions for a
given lattice polygon $P$.

\begin{Th}\label{T:decomp}
Let $P$ be a lattice polygon in $\R^2$ with full Minkowski length $L(P)$.
Consider a maximal Minkowski decomposition in $P$:
$$Q=Q_1+\dots+Q_{L(P)},$$
for some  $Q\subseteq P$.  Then one of the following holds:
\begin{enumerate}
\item  every $Q_i$ is either a primitive segment or a unit triangle;
\item after an  $\A{2}$-transformation and reordering of the summands the
decomposition is
$$Q=T_0+m_1[0,e_1]+m_2[0,e_2]+m_3[0,e_1+e_2],$$
where $m_i$ are non-negative integers such that $m_1+m_2+m_3=L(P)-1$ and the
$e_i$ are the standard basis vectors.
\end{enumerate}
\end{Th}

\begin{pf} Since every $Q_i$ must be strongly indecomposable,  by \rt{indecomp}
it is a primitive segment, a unit triangle, or an exceptional
triangle. We claim that if
one of the $Q_i$ is an exceptional triangle then the other summands are
primitive segments in only three possible directions. This follows from
the two lemmas below.
\end{pf}

\begin{Lemma}\label{L:ml3} Consider two primitive segments $E_1, E_2$
in $\Z^2$ and let $v_1$, $v_2$ be the corresponding vectors. If
$|\det(v_1,v_2)|\geq 3$ then $L(E_1+E_2)\geq 3$.
\end{Lemma}
\begin{pf} We can assume that $v_1=(1,0)$ and $v_2=(a,b)$ with $0\leq
a< b$ and $b=\det(v_1,v_2)$. Cases when $3\leq b\leq 6$ are easily
checked by hand.
For $b\geq 7$ we can use the same argument as in the proof of \rt{indecomp} to
show that $\Pi=E_1+E_2$ contains a segment of lattice length 3.
Indeed, the area of $\Pi$ equals $b\geq 7$. By
Pick's formula $\Pi$ has at least 10 lattice points. But then there exist
$a=(a_1,a_2)$ and $b=(b_1,b_2)$ in~$\Pi$ such that $a_i\equiv b_i$ mod
3, for $i=1,2$.
Therefore the segment  $[a,b]$ is contained in $\Pi$ and has lattice length 3.
\end{pf}

\begin{Lemma}\label{L:triangles}
 Let $P\subset \R^2$ be strongly indecomposable.
Then $L(T_0+P)\geq 3$ unless $P$ is a primitive segment in the direction of
$e_1$,  $e_2$ or $e_1+e_2$.
\end{Lemma}

\begin{pf} Let $E_1$ be an edge of $T_0$ and  $E_2$ an edge of $P$ and
let $v_1,v_2$ be the corresponding vectors. If $|\det(v_1,v_2)|\geq 3$ then by \rl{ml3}
$L(E_1+E_2)\geq 3$ and
since $E_1+E_2\subseteq T_0+P$ we also have $L(T_0+P)\geq 3$.
So we suppose that $|\det(v_1,v_2)|\leq 2$ for all edges $E_1$ in $T_0$.
Then we have the following linear inequalities for $v_2=(s,t)$:
$$
-2\leq s+t\leq 2,\quad -2\leq 2s-t\leq 2,\quad -2\leq s-2t\leq 2.
$$
Clearly, the only integer solutions (up to central symmetry) are
$v_1=(1,0)$, $(0,1)$, and  $(1,1)$.
Now if $P$ contains at least 2 edges in these directions then it must
also contain (up to a translation) either $T=\spn\{(0,0),(1,0),(1,1)\}$
or $T=\spn\{(0,0),(0,1),(1,1)\}$. But in both cases the sum $T_0+T$
contains a $1\times 2$ rectangle which has Minkowski length three.
Therefore, $L(T_0+P)\geq 3$.
\end{pf}

\begin{Rem}\label{R:invTzero} Notice that in \rl{triangles} the
special directions $e_1$,  $e_2$ or $e_1+e_2$ have an easy
$\A2$-invariant description: they are obtained by connecting the interior
lattice point in $T_0$ to the vertices.
\end{Rem}

While classifying polygons of every given full Minkowski length does
not seem feasible, we will make a few statements about polygons of
full Minkowski length 2, which we will use later.


 \begin{figure}[h]
\centerline{
 \scalebox{0.6}
 {
\input{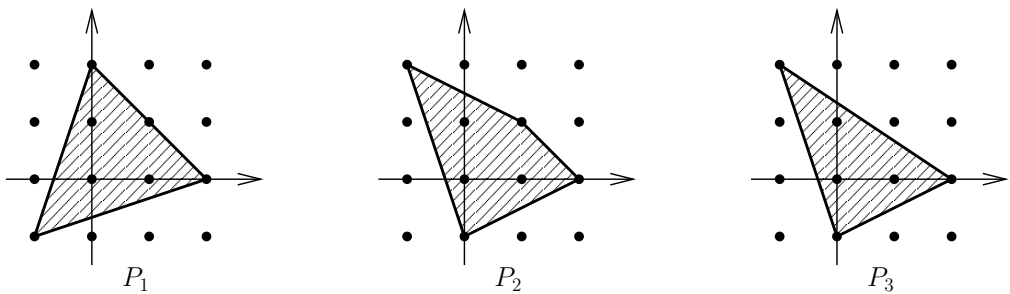}}}
\caption{Full length 2 polygons with 3 interior lattice points.}
\label{F:length2}
\end{figure}

\begin{Prop}\label{P:length2} Suppose $L(P)=2$. Then
\begin{enumerate}
\item $P$ has at most 3 interior lattice points, i.e. $I(P)\leq 3$;
\item if $I(P)=3$ then $P$ is $\A2$-equivalent to one of the polygons
depicted in \rf{length2};
\item if $I(P)=3$ then $L(P+T_0)\geq 4$.
\end{enumerate}
\end{Prop}

\begin{pf} (1) The proof is somewhat technical so we will sketch its
major steps.
Assume $P$ has 4 or more interior lattice points.
First, it is not hard to show that one can choose 4 interior lattice
points in $P$ so that
after an $\A{2}$-transformation they form either a unit square:
$\{(0,0),(1,0),(0,1),(1,1)\}$
or a base 2 isosceles triangle:  $\{(-1,0),(0,0),(1,0),(0,1)\}$.

In the first case, note that $P$ must
include a lattice point which is distance one from the square and lies
on one of the lines containing the sides of the square. By symmetry we can assume it is
$(2,0)$.  In \rf{square} on the left, the solid dots represent the 5
points that now belong to $P$, the crosses represent the points that
cannot belong to $P$  (otherwise its length would be greater than 2).
 \begin{figure}[h]
 \centerline{
 \scalebox{0.52}
 {
\input{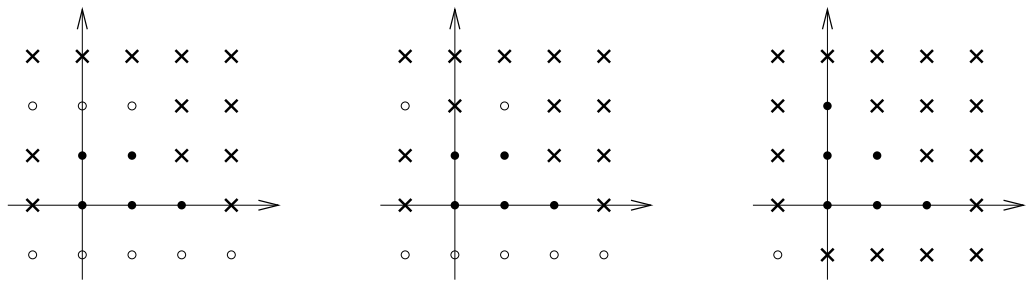}}}
\caption{Non-existence of full length 2 polygons with $I(P)>3$.}
\label{F:square}
\end{figure}
Now if point $(0,2)$ does not belong to $P$ (the middle picture in
\rf{square}) then
either $(-1,2)$ or $(1,2)$ does. But in either case the four points of
the unit square
cannot all lie in the interior of $P$.  If point $(0,2)$ does belong to $P$ then
it produces more forbidden points (the rightmost picture in \rf{square}).  Then
again, it is not hard to see that no such $P$ can exist.

Playing the same game one can show that no $P$ exists in the second
case as well.

(2) First one can show that the three interior lattice points cannot
be collinear. Thus
we can assume that they are $\{(0,0),(1,0),(0,1)\}$.
Our first case is when $(1,1)$ also lies in $P$. Since this must be a
boundary point
and there are no more interior points in $P$ we see that $(-1,2)$ and
$(0,2)$ are
the only possible boundary points of $P$ on the line $y=2$. Similarly,
$(2,0)$ and $(2,-1)$
are the only possible boundary points of $P$ on the line $x=2$. Since
both $(-1,2)$
and $(2,-1)$ cannot belong two $P$, using symmetry we arrive at two
possibilities for the
boundary piece of $P$ containing $(1,1)$, depicted in \rf{3points} on
the left. As in part (1) we crossed
out the points which cannot appear in $P$ since $L(P)=2$.
 \begin{figure}[h]
 \centerline{
 \scalebox{0.52}
 {
\input{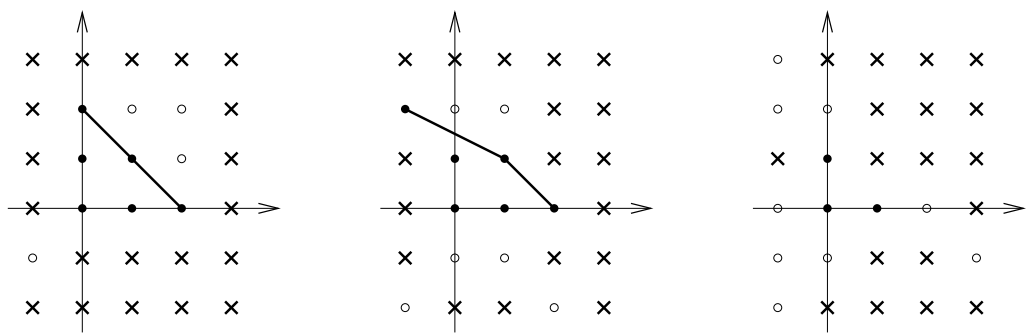}}}
\caption{Constructing full length 2 polygons with $I(P)=3$.}
\label{F:3points}
\end{figure}
Then it becomes clear that
the only $P$ (up to symmetry) containing $\{(0,0),(1,0),(0,1)\}$ and
$(1,1)$ are $P_1$ and $P_2$
in \rf{length2}.

In the second case, when $(1,1)$ does not lie in $P$ we can assume
that $(1,-1)$ and
$(-1,1)$ do not lie in $P$ as well, otherwise we can reduce it to the
previous case by
a unimodular transformation. Also, both $(2,-1)$ and $(-1,2)$ cannot
lie in $P$, so
by symmetry we can assume that $(2,-1)$ does not. As before crossing
out forbidden
points we obtain the rightmost picture in \rf{3points}.
 Now it is easy to see that the only $P$ containing
the 3 points in the interior is $P_3$ in \rf{length2}.

(3) By (2) it is enough to check that $L(P_i+T)\geq 4$ for every
$1\leq i\leq 3$ and
any exceptional triangle $T$.

 We first look at $P_1$.
By \rl{triangles} and \rr{invTzero} we have $L(E+T)\geq 3$ for any
primitive segment $E$ except for the
 three special segments $E_1$, $E_2$, $E_3$ that connect the interior
lattice point of $T$ to its vertices.  If
 $T\neq T_0$ then one of  $[0,e_1]$, $[0,e_2]$, $[0,e_1+e_2]$ is not
among the $E_i$.
 But $P_1$ contains the segments $2[0,e_1], 2[0,e_2]$, and $
(-1,-1)+2[0,e_1+e_2]$.
 If, say, $[0,e_1]$ is not among the $E_i$ then
$L(2[0,e_1]+T)\geq 4$ and hence $L(P_1+T)\geq 4$. It remains to show
that $L(P_1+T_0)\geq 4$
which can easily be checked by hand.

A similar argument works for $P_3$. We only need to replace $T_0$ with
$T_0'$, the triangle with vertices $(0,0), (1,1),$ and $(-1,2)$.
Its special segments  $[0,e_1]$, $[0,e_2]$, $[0,-e_1+e_2]$ are
contained in $P$ with multiplicity 2. Finally, since $P_3\subset P_2$
we do not need to do any extra work for $P_2$.
\end{pf}


\section{Bounds for toric surface codes}\label{S:main}

\subsection{Toric surface codes}
Fix a finite filed $\F_q$ where $q$ is prime power. For any convex
lattice polygon
$P$  in $\R^2$ we associate a $\F_q$-vector space of bivariate polynomials
whose monomials have exponent vectors in $P\cap\Z^2$:
$$\cL(P)=\spn_{\F_q}\{t^m\ |\ m\in P\cap\Z^2\},\quad\text{ where
}t^m=t_1^{m_1}t_2^{m_2}.$$
If  $P$ is contained in the square $K_q^2=[0,q-2]^2$ then the
monomials $t^m$ are
linearly independent over $\F_q$ and so $\dim\cL(P)=|P\cap\Z^2|$. In
what follows we will always
assume that $P\subset K_q^2$.

The {\it toric surface code} $\cC_P$ is a linear code whose codewords
are the strings of values
of $f\in\cL(P)$ at all points of the algebraic torus $(\F_q^*)^2$:
$$\cC_P=\{\left(f(t), t\in(\F_q^*)^2\right)\ |\ f\in\cL(P)\}.$$
This is a linear code of block length $(q-1)^2$ and dimension
$|P\cap\Z^2|$. The
weight of each non-trivial codeword equals the number of points
$t\in(\F_q^*)^2$ where the corresponding polynomial does not vanish.
Let $Z(f)$ denote the number of points in $(\F_q^*)^2$ where $f$ vanishes.
Then the minimum distance $d(\cC_P)$, which is also the minimum
weight, equals
$$d(\cC_P)=(q-1)^2-\max_{0\neq f\in\cL(P)}Z(f).$$

\subsection{The Hasse--Weil bound}\label{S:prelim}
Consider $f\in\cL(P)$. Its {\it Newton polygon} $P_f$ is the convex
hull of the lattice points in $\R^2$ corresponding to the monomials in $f$, so
$$f(t)=\sum_{m\in P_f\cap\Z^2}\lambda_mt^m,\quad\text{where }
t^m=t_1^{m_1}t_2^{m_2}.$$
Let $X$ be a toric variety over $\overline\F_q$ defined by a fan
$\Sig_X$ which is a refinement of the normal fan of $P_f$. Then $f$
can be identified with a global section of a semiample divisor on~$X$.
If $f$ is absolutely irreducible then it defines an irreducible curve
$C_f$ on $X$ whose
number of $\F_q$-rational points $|C_f(\F_q)|$ satisfies the
Hasse--Weil bound (see for example \cite{TVN}):
$$|C_f(\F_q)|\leq q+1+\lfloor 2g\sqrt{q}\rfloor,$$
where $g$ is the arithmetic genus of $C_f$. It is a standard fact from
the theory of toric varieties
that the genus $g$ equals the number $I(P_f)$ of interior lattice
points in $P_f$ (see \cite{Fu}).
Let $D\subset X$ be the invariant divisor at ``infinity'', i.e.
$D=X\setminus (\F_q^*)^2$.
Some of the $F_q$-rational points of $C_f$ may lie on $D$, we will
denote their number by
$B(C_f)$. Then we have the following bound for the number of
$\F_q$-rational points
of $C_f$ in the torus $(\F_q^*)^2$, i.e the number of zeroes of $f$ in
$(\F_q^*)^2$:
\begin{equation}\label{e:HW0}
Z(f)\leq q+1+\lfloor 2g\sqrt{q}\rfloor-B(C_f).
\end{equation}

The divisor $D$ is the disjoint union of zero- and one-dimensional
orbits in $X$.
The one-dimensional orbits $O$ are isomorphic to ${\overline\F_q}^*$ and
correspond to the rays of $\Sig_X$.
Since $\Sig_X$ is a refinement of the normal fan of $P_f$, some of the
orbits correspond to
the edges of $P_f$. Let $E$ be an edge of $P_f$ and $O_E$ the
corresponding orbit in $X$, and consider the ``restriction'' of $f$ to
$E$, i.e. a univariate polynomial $f_E(s)$ whose coefficients
are $\lambda_m$ for $m\in E$, ordered counterclockwise. 
Then the intersection
number $C_f\cdot O_E$ equals the number of zeroes of
$f_E$ in ${\overline\F_q}^*$ (see for example
 \cite{Kho}). In particular if
 $E$ is primitive then $f_E$ is a binomial which has one
$F_q$-rational zero on~$O_E$.
 Therefore, $B(C_f)$ is greater than or equal to the number of
primitive edges of $P_f$. We obtain the
 following proposition.

 \begin{Prop}\label{P:HW}
 Let $f\in\cL(P)$ be absolutely irreducible and $P_f$ its Newton polygon.
 \begin{enumerate}
 \item If $P_f$ is an exceptional triangle then $Z(f)\leq q-2+\lfloor
2\sqrt{q}\rfloor$.
 \item If $I(P_f)=0$ then $Z(f)\leq q-1$ unless $P_f$
 is twice a unit triangle in which case $Z(f)\leq q+1$.
 \end{enumerate}
\end{Prop}

\begin{pf} (1) follows immediately from \re{HW0} and the above
discussion. For (2) we use the
classification of polygons with no interior lattice points (see for
example \cite{BaNi}):
$P_f$ is $\A{2}$-equivalent to either (a) $2\D$ or (b) a trapezoid
(see \rf{empty})
where $0\leq a\leq b$ (this includes primitive segments when $a=b=0$
and unit triangles when
$a=0$, $b=1$).
 \begin{figure}[h]
\centerline{
 \scalebox{0.6}
 {
\input{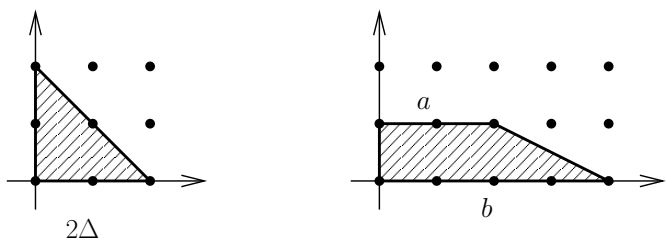}}}
\caption{Polygons with no interior lattice points.}
\label{F:empty}
\end{figure}
In the first case $Z(f)\leq q+1$ by \re{HW0}.
In the second case $P_f$ has at least 2 primitive edges, so $Z(f)\leq
q-1$, again by \re{HW0}.
\end{pf}

\subsection{Bounds for the minimum distance}
Let $\cC_P$  be the toric surface code defined by a lattice
polygon~$P$ in~$K_{q}^2$.
In this section we
prove bounds for the minimum distance of $\cC_P$ in terms of the
full Minkowski length $L(P)$
of the polygon $P$.

Here is our first application of the results of the previous section.

\begin{Prop}\label{P:HWL} 
Let  $f\in\cL(P)$ be a polynomial with the largest number of
absolutely irreducible factors, $f=f_1\cdots f_L$.
Then
\begin{enumerate}
\item $L=L(P)$ and every $P(f_i)$ is either a primitive segment, a
unit triangle, or an exceptional triangle;
\item the number of zeroes of $f$ in $(\F_q^*)^2$ satisfies
$$Z(f)\leq L(q-1)+\lfloor 2\sqrt{q}\rfloor-1;$$
\item if $P(f_i)$ is not an exceptional triangle for any $1\leq i\leq L$
then $$Z(f)\leq L(q-1).$$
\end{enumerate}
\end{Prop}

\begin{pf} (1) follows directly from \rt{decomp}. Moreover, the theorem implies that
either (a) all $P_i$ are primitive segments
or unit triangles
or (b) one of the $P_i$ is an exceptional triangle and the others are
primitive segments.

 In the first case every $f_i$ has at most $q-1$ zeroes in
$(\F_q^*)^2$ by \rp{HW}.
 Not accounting for possible common zeroes of the $f_i$ we obtain the
bound in (3).
 In the second case one of the $f_i$ has at most $q-2+\lfloor
2\sqrt{q}\rfloor$ zeroes
 and the others have at most $q-1$ zeroes, again by \rp{HW}.
 As before, disregarding
possible common zeroes of the $f_i$ we get the bound in (2).
\end{pf}

The next proposition deals with polynomials $f$ whose number of
absolutely irreducible
factors is $L(P)-1$.

\begin{Prop}\label{P:HWL-1}
Let $P$ have full Minkowski length $L$ and let $f\in\cL(P)$ have $L-1$
absolutely irreducible factors. Then
$$Z(f)\leq (L-1)(q-1)+\lfloor 6\sqrt{q}\rfloor.$$
\end{Prop}

\begin{pf}
As before let $f=f_1\cdots f_{L-1}$ be the decomposition of $f$ into absolutely
irreducible factors and let $P_i$ be the Newton polygon of $f_i$.
First, by \rp{simple}
$$k+1=L\geq\sum_{i=1}^{k} L(P_i)\geq k,$$
hence,
up to renumbering, $L(P_1)\leq 2$ and $L(P_i)=1$ for $2\leq i\leq k$.

Assume $L(P_1)=1$. Then
every $P_i$ is either a strongly indecomposable triangle or a lattice
segment. We claim that at most 3 of the $P_i$ are exceptional
triangles, and so the statement
follows from \rp{HW}. Indeed, if, say,  $P_1,\dots, P_4$ are
exceptional triangles
 then by \rl{triangles} $L(P_1+\dots+P_4)\geq 6$. Applying
\rp{simple} again we get
$$k+1=L \geq L(P_1+\dots+P_4)+\sum_{i=5}^kL(P_i)\geq 6+(k-4)=k+2,$$ a
contradiction.

Now assume $L(P_1)=2$. According to \rp{length2}, (1) we have $I(P_1)\leq 3$.
Also since $L(P_1)=2$, at most one of the other $P_i$ is an
exceptional triangle.
This follows from \rl{triangles} using similar to the previous case arguments.
We now have three subcases.

If $I(P_1)=1$ then we have
$$Z(f)\leq (q+1+\lfloor 2\sqrt{q}\rfloor)+(q-2+\lfloor
2\sqrt{q}\rfloor)+(L-3)(q-1)\leq (L-1)(q-1)+\lfloor
6\sqrt{q}\rfloor.$$

If $I(P_1)=2$ then $P_1$ has at least one primitive edge which we
prove in the lemma below.
Therefore by \re{HW0} and \rp{HW} we have
$$Z(f)\leq (q+\lfloor 4\sqrt{q}\rfloor)+(q-2+\lfloor
2\sqrt{q}\rfloor)+(L-3)(q-1)\leq (L-1)(q-1)+\lfloor
6\sqrt{q}\rfloor.$$

Finally, if $I(P_1)=3$ then none of the other $P_i$ is an exceptional
triangle. This follows from
\rp{length2}, (3) and the above arguments. In this case $P_1$ has at
least 2 primitive edges
by \rp{length2}, (2). Therefore using \re{HW0} we have
$$Z(f)\leq (q-1+\lfloor 6\sqrt{q}\rfloor)+(L-2)(q-1)=
(L-1)(q-1)+\lfloor 6\sqrt{q}\rfloor.$$

\begin{Lemma} If $L(P)=2$ and $I(P)=2$ then $P$ has a primitive edge.
\end{Lemma}
\begin{pf} Since $L(P)=2$ no edge can have more than 3 lattice points.
If $P$ has 4 or more edges none of which is primitive then $P$
has at least 8 boundary lattice points and, hence, at least 10 lattice points total. But then $P$
contains a lattice segment of lattice length 3 (see the proof of \rl{ml3}), which contradicts the
assumption $L(P)=2$.

It remains to show that triangles with no primitive edges, 2 interior
lattice points, and 6 boundary lattice points do not exist. Let $T$ be such a triangle
and let $2E_1$, $2E_2$ be two of its edges, where $E_1$ and $E_2$
are primitive. Then $E_1$, $E_2$ form a triangle $T'$
of area $A(T')=\frac{1}{4}A(T)$. On the other hand, by Pick's formula $A(P)=4$, and hence
$A(T')=1$. This implies that up to an $\A{2}$-transformation $E_1=[0,e_1]$ and $E_2=[0,e_1+2e_2]$,
but then $I(T)=1$, a contradiction.
\end{pf}

\end{pf}

Now we are ready for the main result of this section.

\begin{Th}\label{T:2-bound}
Let $P\subset K_{q-1}^2$ be a lattice polygon with area $A=A(P)$ and
full Minkowski length $L=L(P)$.
Then
\begin{enumerate}
\item for $q\geq
\max\left(23,\big(c+\sqrt{c^2+5/2}\big)^2\right)$,
where $c=A/2-L+9/4$,
every polynomial $f\in\cL(P)$ has  at most $L(q-1)+2\sqrt{q}-1$ zeroes in
$(\F_q^*)^2$.
Consequently, the minimum distance for the toric surface code $\cC_P$ satisfies
$$d(\cC_P)\geq (q-1)^2-L(q-1)-2\sqrt{q}+1.$$
\item if no maximal decomposition in $P$ contains an exceptional
triangle then for
$q\geq \max\left(37,\big(c+\sqrt{c^2+2}\big)^2\right)$,
where $c=A/2-L+11/4$,
every polynomial $f\in\cL(P)$ has at most
$L(q-1)$ zeroes in $(\F_q^*)^2$.
Consequently, the minimum distance for the toric surface code $\cC_P$ satisfies
$$d(\cC_P)\geq (q-1)^2-L(q-1).$$

\end{enumerate}
\end{Th}

\begin{pf}
(1) As we have seen in \rp{HWL} (2) the bound holds for
the polynomials with the largest number of irreducible factors. We are
going to show that
for large enough~$q$ every polynomial with fewer irreducible factors
will have no greater than $L(q-1)+\lfloor 2\sqrt{q}\rfloor-1$ zeroes
in $(\F_q^*)^2$.

Let $f\in\cL(P)$ have $k<L$ absolutely irreducible factors
$f=f_1\cdots f_k$ and let $P_i$
be the Newton polygon of $f_i$. If $k=L-1$ we can
use the bound in \rp{HWL-1}:
\begin{equation}\label{e:ourbound1}
Z(f)\leq (L-1)(q-1)+\lfloor 6\sqrt{q}\rfloor.
\end{equation}
The latter is at most $L(q-1)+\lfloor 2\sqrt{q}\rfloor-1$ for all $q\geq 19$.

Now suppose $1\leq k\leq L-2$. First assume $I(P_i)=0$ for all $1\leq i\leq k$.
Then by \rp{HW} (2),
$$Z(f)\leq s(q+1)+(k-s)(q-1)=2s+k(q-1),$$
where $s$ is the number of twice unit triangles among the $P_i$. Since the sum
of the full Minkowski lengths of the $P_i$ cannot exceed $L$ we have
$2s+(k-s)\leq L$, i.e. $s\leq L-k$.
Using this inequality along with $k\leq L-2$ we obtain
$$Z(f)\leq 2s+k(q-1)\leq 2L+k(q-3)\leq (L-2)(q-1)+4.$$
The latter is at most $L(q-1)$ for all $q\geq 3$ and the bounds follow.

Suppose $I(P_i)>0$ for at least one of the $P_i$. Then, as we will
show in \rl{others},
\begin{equation}\label{e:ourbound2}
Z(f)\leq k(q-1)+2\big(A+3/2-2k\big)\sqrt{q}+2.
\end{equation}
Now the right hand side will be at most $L(q-1)+2\sqrt{q}-1$ whenever
$q$ satisfies
\begin{equation}\label{e:quadratic}
(L-k)q-2(A+1/2-2k)\sqrt{q}-(L-k+3)\geq 0.
\end{equation}
Before proceeding we introduce the following notation: $m=L-k$,
$d=A/2-L+1/4$. Then \re{quadratic} becomes
$$mq-4(d+m)\sqrt{q}-(m+3)\geq 0,\quad 2\leq m\leq L-1.$$
Since this is a quadratic inequality in $\sqrt{q}$, it will hold if
$$\sqrt{q}\geq C+\sqrt{C^2+1+3/m},\quad \text{where }C=2+2d/m.$$
Since $m\geq 2$ it is enough to choose $\sqrt{q}\geq C+\sqrt{C^2+5/2}$.
Finally, if $d\geq 0$ then $C\leq 2+d$, since $m\geq 2$, and it is
enough to choose
$$q\geq \big(c+\sqrt{c^2+5/2}\big)^2,\quad\text{where}\quad c=2+d=A/2-L+9/4.$$
If $d<0$ then $C< 2$ and it is enough to choose $q\geq 23$.

(2) The proof of the second statement is completely analogous. First,
if $f$ has $L$ irreducible factors the bound holds by \rp{HWL}, (3).
Second, if $f$ has fewer than $L$
factors we choose $q$ large enough so that the right hand sides of
\re{ourbound1} and \re{ourbound2}
are no greater than $L(q-1)$. The same arguments as before show that
it is enough to
choose
$$q\geq \max\left(37,\big(c+\sqrt{c^2+2}\big)^2\right),\quad\text{where}\quad
c=A/2-L+11/4.$$
\end{pf}

It remains to prove the following lemma.

\begin{Lemma}\label{L:others} Let $f=f_1\cdots f_k$,  for $1\leq k\leq
L-2$, and $I(P_i)>0$
for at least one $i$. Then
$$Z(f)\leq k(q-1)+2\big(A+3/2-2k\big)\sqrt{q}+2.$$
\end{Lemma}
\begin{pf}
We order the $P_i$ so that for $1\leq i\leq t$ every $P_i$
either has interior lattice points or is twice a unit triangle.
 Then according to \re{HW0}
and \rp{HW} we have
\begin{equation}\label{e:before}
Z(f)\leq t(q+1)+2\sqrt{q}\sum_{i=1}^tI(P_i)+(k-t)(q-1).
\end{equation}
Now we want to get a bound for $\sum_{i=1}^tI(P_i)$. Recall that given
two polytopes $Q_1$ and $Q_2$ in $\R^2$, their normalized {\it mixed
volume} (2-dimensional) is
$$V(Q_1,Q_2)=A(Q_1+Q_2)-A(Q_1)-A(Q_2).$$
The mixed volume is symmetric; bilinear with respect to Minkowski
addition; monotone increasing
(i.e. if $Q_1'\subset Q_1$ then $V(Q_1',Q_2)\leq V(Q_1,Q_2)$); and
$\A{2}$-invariant (see, for example
\cite{BuZa}, page 138).
This implies that
\begin{equation}\label{e:mixed}
V(P_i,P_j)\geq 2,\quad \text{for }\ 1\leq i\leq t\text{ and }1\leq j\leq k.
\end{equation}
Indeed, by monotonicity it is enough to show that $V(P_i,E)\geq 2$ for
any lattice segment $E$,
and by $\A2$-invariance we can assume that $E$ is horizontal. It
follows readily from the definition
that  $V(P_i,E)=h(P_i)|E|$, where $h(P_i)$ is the length of the
horizontal projection of $P_i$
(the height of $P_i$) and $|E|$ is the length of~$E$. Clearly,
$|E|\geq 1$  and $h(P_i)\geq 2$ if $P_i$ has at least one interior
lattice point or is twice a unit triangle.

Using \re{mixed} and bilinearity of the mixed volume, by induction we obtain
\begin{eqnarray}
A&\geq &A\Big(\sum_{i=1}^kP_i\Big)=A(P_1)+A
\Big(\sum_{i=2}^kP_i\Big)+V \Big(P_1,\sum_{i=2}^kP_i \Big)\nonumber\\
&\geq&A(P_1)+A \Big(\sum_{i=2}^kP_i \Big)+2(k-1)\geq\dots\nonumber\\
&\geq&\sum_{i=1}^tA(P_i)+A
\Big(\sum_{i=t+1}^kP_i\Big)+2\sum_{i=1}^t(k-i)\geq\sum_{i=1}^tA(P_i)+2kt-t^2-t.\nonumber
\end{eqnarray}
Now, by Pick's formula $A(P_i)=I(P_i)+\frac{1}{2}B(P_i)-1\geq
I(P_i)+\frac{1}{2}$ since $B(P_i)$,
the number of boundary lattice points, is at least 3. Therefore
$$\sum_{i=1}^tI(P_i)\leq A+t^2+\frac{t}{2}-2kt.$$
Substituting this into \re{before} and simplifying we obtain
\begin{equation}\label{e:after}
Z(f)\leq k(q-1)+2\sqrt{q}\big(A+t^2+\frac{t}{2}-2kt\big)+2t.
\end{equation}
It remains to note that the maximum of the right hand side of \re{after}
is attained at $t=1$, provided $k\geq 1$ and $q\geq 4$, and that establishes the required inequality.
\end{pf}


\section{Two algorithms}\label{S:alg}

Given a polytope $P$, to make use of our bound in \rt{2-bound} it
remains to understand
\begin{enumerate}
\item how to find $L(P)$, the full Minkowski length of $P$,
\item how to determine  whether there is a maximal Minkowski
decomposition in $P$ one of whose summands is an exceptional triangle.
\end{enumerate}

Here we provide algorithms that answer these questions in polynomial
time in $|P\cap\Z^2|$.

Recall that a zonotope  $Z=\sum_{i=1}^k E_j\subseteq P$ is called {\it
maximal for $P$} if  $k$, the number of non-trivial  Minkowski
summands  (counting
their multiplicities), is equal to $L(P)$.

It follows from \rp{simple} that a maximal zonotope always exists
although it is usually not unique. It turns
out that any maximal zonotope of $P$ has at most four distinct
summands  and among them there are maximal zonotopes with a
particularly  easy description.

\begin{Prop}\label{P:3segments}
Let $P$ be a lattice polygon. Then
\begin{enumerate}
\item Any zonotope $Z$ maximal for $P$ has  at most 4 different summands.
\item There exists a zonotope $Z$ maximal for $P$ with at most 3
different summands. Moreover up to an $\A{2}$-transformation
these summands are $[0,e_1]$, $[0,e_2]$, and $[0,e_1+e_2]$.
\end{enumerate}
\end{Prop}

\begin{pf} Let $Z=\sum_{i=1}^{L} E_j$  be a  zonotope maximal for
$P$ and let $v_j$ be the vector of $E_j$. According to \rl{ml3},
$|\det(v_i,v_j)|\leq 2$ for any $1\leq i, j\leq k$.

The case when all $v_i$ are the same is trivial. Suppose there are
exactly two different summands,
i.e. $Z=m_1E_1+m_2E_2$ for some positive integers $m_1\geq m_2$ and
$E_1\neq E_2$.
If $|\det(v_1,v_2)|=1$ then we
can transform $(v_1,v_2)$ to the standard basis $(e_1,e_2)$ and (2) follows.
If $|\det(v_1,v_2)|=2$ then we
can assume that $v_1=e_1$ and  $v_2=e_1+2e_2$. But $E_1+E_2$ contains
$2[0,e_2]$, so we can
pass to $Z'=(m_1-m_2)[0,e_1]+2m_2[0,e_2]$. Clearly $Z'\subseteq Z$
and $Z'$ is maximal.

Now suppose that $Z$ has at least three different summands.
First, assume $|\det(v_i,v_j)|=2$ for all $i\neq j$.  As before,
without loss of generality,
$v_1=e_1$ and  $v_2=e_1+2e_2$. Consider $v_3=(s,t)$. By looking at the
determinants $\det(v_i,v_3)$
for $i=1,2$ we have $|t|=2$ and $|t-2s|=2$. This implies that $v_3$ is
not primitive, a contradiction. Therefore, $|\det(v_i,v_j)|=1$ for
some $i\neq j$ and
we can assume that $v_1=e_1$ and  $v_2=e_2$. Again, we let
$v_3=(s,t)$ and look at the
determinants $\det(v_i,v_3)$ for $i=1,2$.  We see that the only
vectors $v_3$ (up to central symmetry) that may appear are $(1,1)$,
$(1,-1)$, $(2,1)$, $(2,-1)$, $(1,2)$, $(1,-2)$. No two out of the last
four vectors can appear together as they generate parallelograms of
area at least ~3. For the same reason $(1,1)$ cannot appear with
$(2,-1)$ or $(1,-2)$, and $(1,-1)$ cannot appear with $(2,1)$ or
$(1,2)$. We have three possible combinations:
\begin{itemize}
\item[(a)] $v_1=(1,0), v_2=(0,1), v_3=(1,1), v_4=(1,-1)$
\item[(b)] $v_1=(1,0), v_2=(0,1), v_3=(1,1)$ and $v_4=(1,2)$  or $v_4=(2,1)$
\item[(c)] $v_1=(1,0), v_2=(0,1), v_3=(1,-1)$ and $v_4=(1,-2)$  or $v_4=(2,-1)$
\end{itemize}

We have proved our first claim. To prove the second, note that we
can actually reduce the number of distinct segments $E_j$. In case (a)
$2E_1\subset E_3+E_4$ and we will be able to get rid of either $E_3$
or $E_4$ by replacing $E_3+E_4$ with $2E_1$. In either case the remaining
segments  are $\A2$-equivalent to $[0,e_1]$, $[0,e_2]$, and
$[0,e_1+e_2]$.

In case (b) we can assume that $v_4=(1,2)$.  Since $2E_2\subset
E_1+E_4$ we will be able to get rid of either $E_1$ or $E_4$ and the
remaining segments  are $\A2$-equivalent to $[0,e_1]$, $[0,e_2]$, and
$[0,e_1+e_2]$.  Case (c) is obtained from (b) by flipping the second
coordinate.
\end{pf}

To find $L(P)$ we only need to look at all the zonotopes  $Z\subseteq
P$ with at most three different summands $\A{2}$-equivalent to
$[0,e_1]$, $[0,e_2]$, and $[0,e_1+e_2]$ and find the one that has the
largest number of summands (counting multiplicities).

\begin{Th}\label{T:alg1}
 Let $P$ be a lattice polygon and  let $|P\cap\Z^2|$ be the number of
lattice points in $P$. Then the full Minkowski length $L(P)$ can be
found in polynomial time in $|P\cap\Z^2|$.
\end{Th}
\begin{pf}
The case when $P$ is 1-dimensional is trivial so we will be assuming
that $P$ has dimension 2.

For every triple of points $\{A,B,C\}\subseteq P\cap\Z^2$, where it is
important which point goes first and the order of the other two does
not matter, we check if $E_1=[A,B]$ and $E_2=[A,C]$ generate a
parallelogram of area one. If so,  we want to construct various
zonotopes  whose summands are $E_1$, $E_2$ and $E_3=[A,B+C]$. We do
this in the most straightforward way.

First, for every $1\leq i\leq 3$, we find $M_i$, the largest integer
such that a lattice translate of  $M_iE_i$ is contained in $P$.  For
this we find the maximum number of lattice points in the linear
sections of $P$
with lines in the direction of $E_i$ (there are finitely many such
lines with at least one lattice point of $P$).

Second, for each triple of integers $m=(m_1, m_2, m_3)$ where $0\leq
m_i\leq M_i$,
we check if some lattice translate of the zonotope
$Z_m=m_1E_1+m_2E_2+m_3E_3$ is contained in $P$ (we run through lattice
points $D$ in $P$ to check if  $D+Z_m$ is contained in $P$). For all
such zonotopes that fit into $P$ we look at $m_1+m_2+m_3$ and find the
maximal possible value $M$ of this sum.

Finally, the largest such sum $M$ over all choices of
$\{A,B,C\}\subseteq P\cap\Z^2$ is $L(P)$, by
\rp{3segments}. Clearly, this algorithm is polynomial in $|P\cap\Z^2|$.

Notice that in the above we have taken care of the maximal zonotopes
that are possibly multiples of a single segment.  Indeed, if $[A,B]$
is a primitive segment connecting two lattice points in $P$ then
unless $P$ is 1-dimensional  there is a lattice point $C$ in $P$ such
that $[A,B]$ and $[A,C]$ generate a parallelogram of area one.  We can
assume that $A$ is the origin and $B=(1,0)$. Let $C=(k,l)$ be a
lattice point in $P$ with smallest positive $l$ (flip $P$ with respect
to the $x$-axis if necessary).  By the minimality of $l$ the triangle
$ABC$ has no lattice points except its vertices. By Pick's formula,
its area is $1/2$ and we have found the required third vertex $C$.
\end{pf}

\begin{Th}\label{T:alg2}
Let $P$ be a lattice polygon in $\R^n$.
Then we can decide in polynomial time in $|P\cap\Z^2|$ if there is a
maximal Minkowski decomposition in $P$ one of  whose summands  is
an exceptional triangle.
\end{Th}
\begin{pf}
We first run the algorithm from \rt{alg1} to find $L(P)$. Next for
each triple of points $A,B,C\in P\cap\Z^2$ we check if the triangle
$T_{ABC}$ has exactly four lattice points --- the three vertices $A,B,C$
and one point  $D$ strictly inside the triangle. If so, this triangle
is exceptional. If this triangle is a summand in some
maximal Minkowski decomposition in $P$ then the other summands that
may appear in this decomposition are the primitive segments $E_1$, $E_2$,
and $E_3$ connecting $D$ to the vertices $A, B, C$ (see \rr{invTzero}).

Now it remains to look at all Minkowski  sums
$T_{ABC}+m_1E_1+m_2E_2+m_3E_3$ with $m_1+m_2+m_3=L(P)-1$ and check if
any of them fits into $P$. If this indeed happens for some $T_{ABC}$,
there is a maximal decomposition in $P$  with an exceptional triangle.
Otherwise any maximal decomposition is a
sum of primitive segments and unit triangles. Clearly, this algorithm
is polynomial in
$|P\cap\Z^2|$.


\section{Three Examples}\label{S:ex}

In this section we illustrate our methods with three examples. \rex{2} was
given by Joyner in \cite{Jo}. \rex{3} appears in the Little and
Schenck's paper \cite{LiSche}.

\begin{Ex}\label{ex:1} Consider the pentagon $P$ with vertices
$(0,0)$, $(1,0)$, $(1,3)$, $(2,4)$, and $(4,2)$ as in \rf{ex1}.
 \begin{figure}[h]
\centerline{
 \scalebox{0.6}
 {
\input{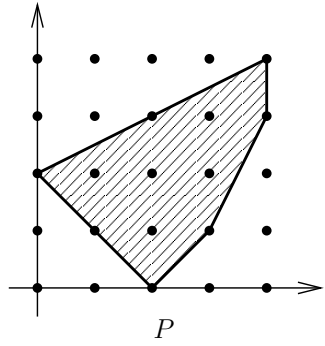}}}
\caption{Pentagon}
\label{F:ex1}
\end{figure}
One can easily check that $L(P)=3$ and there is a maximal
decomposition in $P$ containing
$T_0$. In fact, $P$ contains $T_0+[0,e_2]+[0,e_1+e_2]$. It defines a
toric surface code
of dimension $n=|P\cap\Z^2|=12$. To apply \rt{2-bound}
we compute $A=15/2$, so $c=3$. Therefore,
$$d(C_P)\geq (q-1)^2-3(q-1)-2\sqrt{q}+1,$$
for all $q\geq 41$.
In this particular example we can establish a better lower bound
for~$q$, namely $q\geq 19$.
Indeed, we have already seen in the proof of \rt{2-bound}
that every $f$ with 2 absolutely irreducible factors will have at
 most $3(q-1)+2\sqrt{q}-1$ for all $q\geq 19$ (see \re{ourbound1}).
If $f$ is absolutely irreducible we use \re{HW0}. Then it has at most
$q+1+\lfloor 10\sqrt{q}\rfloor-2$ zeroes since $P_f\subseteq P$ has at most 5 interior
lattice points in which case it will have at least 2 primitive edges. But
$$q+1+\lfloor 10\sqrt{q}\rfloor-2\leq 3(q-1)+2\sqrt{q}-1$$
for all $q\geq 19$.
\end{Ex}

\begin{Ex}\label{ex:2} Consider the triangle $P$ with vertices
$(0,0)$, $(4,1)$, and $(1,4)$ (see \rf{ex2}).
 \begin{figure}[h]
\centerline{
 \scalebox{0.6}
 {
\input{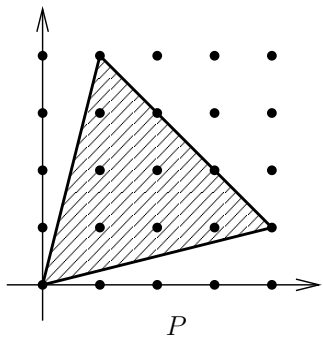}}}
\caption{Triangle}
\label{F:ex2}
\end{figure}
This example is similar to the
previous one. We also have $L(P)=3$,  $A=15/2$, but the dimension of
the corresponding toric surface code is slightly smaller,
$n=|P\cap\Z^2|=11$. However in this case $P$
has no exceptional triangles in any maximal decomposition. Therefore,
\rt{2-bound} provides a
better bound for the minimum distance:
$$d(C_P)\geq (q-1)^2-3(q-1),$$
which holds for all $q\geq 53$. As before, this can be improved to
$q\geq 37$ using \re{HW0} and the fact that
$I(P)=6$. Note that $f=xy(x-a)(x-b)(x-c)$, for $a,b,c\in\F_q^*$ distinct,
has exactly $3(q-1)$ zeroes in  $(\F_q^*)^2$, hence
for $q\geq 37$ the above bound is exact:
\begin{equation}\label{e:exact}
d(C_P)=(q-1)^2-3(q-1).
\end{equation}
For $q=8$ this was previously established by Joyner \cite{Jo}.
Also \re{exact} follows from Little and Schenck's result \cite{LiSche}  for all 
$q\geq (4I(P)+3)^2=729$.
\end{Ex}

\begin{Ex}\label{ex:3}
Let $P$ be the hexagon with vertices $(1,0)$, $(0,1)$, $(1,2)$,
$(3,3)$, $(3,2)$,
and $(2,0)$ (see \rf{ex3}).
 \begin{figure}[h]
\centerline{
 \scalebox{0.6}
 {
\input{ex3.pstex_t}}}
\caption{Hexagon}
\label{F:ex3}
\end{figure}
We have $L(P)=3$, $A=5$ and $\cC_P$ has dimension $9$. Also $P$ has no
maximal decomposition with an exceptional triangle. Therefore,
\rt{2-bound} implies
$$d(C_P)\geq (q-1)^2-3(q-1),$$
for all $q\geq 37$. Little and Schenck's result  \cite{LiSche} proves this bound for $q>225$.
In fact we can show more in this example: for all $q\geq 11$
\begin{equation}\label{e:more}
d(C_P)=(q-1)^2-3(q-1)+2.
\end{equation}
To see this, first note that $f=x(x-a)(y-b)(y-c)$, for
$a,b,c\in\F_q^*$ distinct, has exactly
$3(q-1)-2$ zeroes in $(\F_q^*)^2$. Furthermore, every maximal
decomposition in $P$
is of the form $E_1+2E_2$, where $E_i$ is a primitive segment in the
direction of
$e_1$, $e_2$, or $e_1+e_2$. This implies that every polynomial $f$
with the largest
number of absolutely irreducible factors (three) will have at most
$3(q-1)-2$ zeroes in  $(\F_q^*)^2$
(here we take into account the intersections of the irreducible curves
defined by the factors of $f$).

Now we claim that for $q\geq 11$ polynomials with fewer factors (one
or two) will have at most $3(q-1)-2$ zeroes in  $(\F_q^*)^2$ as well. Indeed, decompositions
with 2 summands in $P$ can have at most one exceptional triangle, hence, 
$Z(f)\leq 2(q-1)+\lfloor2\sqrt{q}\rfloor$ for every $f$ with 2 irreducible factors.
This will be no greater than $3(q-1)-2$ for $q\geq 9$.  If $f$ is absolutely
irreducible then by \re{HW0} $Z(f)\leq q+1+\lfloor 6\sqrt{q}\rfloor-3$, which is
no greater than $3(q-1)-2$ starting with $q=11$.

The computations  preformed in \cite{LiSche} show the validity of \re{more} for all $5\leq q\leq 11$ except for $q=8$ when  the answer is $d(C_P)=(q-1)^2-3(q-1)$. We now have a
complete understanding of this example.
\end{Ex}

\end{pf}



\begin{thebibliography}{9}

\bibitem[1]{BaNi} V. Batyrev, B. Nill, {\em Multiples of lattice polytopes without interior lattice points}, Mosc. Math. J.  7  (2007),  no. 2, 195--207, 349.

\bibitem[2]{BeRo} M. Beck, S. Robins, {\em Computing the continuous discretely. Integer-point enumeration in polyhedra.} Undergraduate Texts in Mathematics. Springer, New York, 2007.

\bibitem[3]{BuZa} Yu. D. Burago, V. A. Zalgaller {\em Geometric
Inequalities}, Springer--Verlag, Berlin, 1988

\bibitem[4]{Fu} W. Fulton, {\em Introduction to Toric Varieties},
Princeton University Press, Princeton, NJ, 1993.

\bibitem[5]{Han} J. Hansen, {\em Toric Surfaces and Error--correcting
    Codes} in Coding Theory, Cryptography, and Related Areas, Springer
  (2000), pp. 132-142.

\bibitem[6]{Jo} D. Joyner, {\em Toric codes over finite fields}, Appl.
Algebra Engrg. Comm. Comput.,
15 (2004), pp. 63--79.

\bibitem[7]{Kho} A. G. Khovanskii, {\em Newton polytopes, curves on toric surfaces, and inversion of Weil's theorem}, Russian Math. Surveys  52  (1997),  no. 6, 1251--1279.

\bibitem[8]{LiSche} J. Little, H. Schenck,
{\em Toric Surface Codes and Minkowski sums},
SIAM J. Discrete Math.  20  (2006),  no. 4, 999--1014 (electronic).

\bibitem[9]{TVN} M. Tsfasman, S. Vl\u adu\c t, D. Nogin,
{\em Algebraic geometric codes: basic notions},
Mathematical Surveys and Monographs, 139,  AMS, Providence, RI, 2007.



\end{thebibliography}
\end{document}